\documentclass[11pt,reqno]{amsart}

\usepackage{amsmath, amsfonts, amsthm, amssymb}

\newtheorem{Theorem}{Theorem}[section]
\newtheorem{Lemma}{Lemma}[section]
\newtheorem{Proposition}{Proposition}[section]

\theoremstyle{definition}

\newtheorem{Remark}{Remark}[section]

\numberwithin{equation}{section}

\renewcommand{\u}{{\bf u}}

\newcommand{\R}{{\mathbb R}}
\newcommand{\Dv}{{\rm div}}

\newcommand{\m}{{\bf m}}

\def\f{\frac}
\renewcommand{\O}{\Omega}

\def\D{\Delta }

\def\hf1{^\f{1}{1-\xi^2}}

\def\be{\begin{equation}}
\def\en{\end{equation}}
\def\bs{\begin{split}}
\def\es{\end{split}}

\renewcommand{\d}{{\bf d}}

\author{Dehua Wang}
\address{Department of Mathematics, University of Pittsburgh,
                           Pittsburgh, PA 15260.}
\email{dwang@math.pitt.edu}
\author{Cheng Yu}
\address{Department of Mathematics, University of Pittsburgh,
                           Pittsburgh, PA 15260.}
\email{chy39@pitt.edu}

\title[incompressible limit  of liquid crystals]
{Incompressible limit  for  the compressible flow of liquid crystals}


\date{\today}

\begin{document}

\begin{abstract}
The connection between the  compressible flow of liquid crystals with low Mach number and the incompressible
flow of liquid crystals is studied in a bounded domain.
In particular, the convergence of weak solutions of the compressible flow of liquid crystals to the weak solutions of the incompressible flow of liquid crystals is proved when the Mach number approaches  zero;
that is, the  incompressible limit is justified for weak solutions in a bounded domain.
\end{abstract}

\maketitle

\section{Introduction}

In this paper, we consider the incompressible limit of the following hydrodynamic system of partial differential equations for the
three-dimensional compressible flow of nematic liquid crystals \cite{DE,HK,Lin2}:
\begin{subequations}\label{I1}
\begin{align}
&\widetilde{\rho_{t}}+\Dv(\widetilde{\rho}\widetilde{\u})=0,\label{I1a}\\
&(\widetilde{\rho}\widetilde{\u})_{t}+\Dv(\widetilde{\rho}\widetilde{ \u} \otimes \widetilde{\u})+\nabla \widetilde{P}(\widetilde{\rho})= \widetilde{\mu} \Delta\widetilde{\u}-\widetilde{\lambda }\Dv\left(\nabla \widetilde{\d} \odot \nabla \widetilde{\d}-\big(\frac{1}{2} |\nabla\widetilde{\d}|^{2}+F(\widetilde{\d})\big) I_{3}\right),\label{I1b}\\
&\widetilde{\d_{t}}+\widetilde{\u}\cdot\nabla\widetilde{\d}=\widetilde{\theta}(\Delta\widetilde{\d}-f(\widetilde{\d})),\label{I1c}
\end{align}
\end{subequations}
where  $\widetilde{\rho}\ge 0$ denotes the density, $\widetilde{\u}\in\R^3$  the velocity, $\widetilde{\d}\in\R^3$ the direction
field for the averaged macroscopic molecular orientations, and
$\widetilde{P}=a\widetilde{\rho}^{\gamma}$ is the pressure with constants $a>0$ and $\gamma\geq 1$.  The positive constants $\widetilde{\mu}, \widetilde{\lambda}, \widetilde{\theta}$ denote  the viscosity, the competition between kinetic energy and potential energy, and the microscopic elastic relation time for the molecular orientation field, respectively.  The symbol $\otimes$ denotes the Kronecker tensor product,
$I_{3}$ is the $3\times3$ identity matrix, and $\nabla\widetilde{\d}\odot\nabla\widetilde{\d}$ denotes the $3\times 3$ matrix whose $ij$-th entry is $<\partial_{x_i} \widetilde{\d}, \partial_{x_j}\widetilde{\d}>$.  Indeed, $$\nabla \widetilde{\d} \odot\nabla \widetilde{\d}= (\nabla\widetilde{\d})^{\top}\nabla\widetilde{\d},$$
where $(\nabla\widetilde{\d})^{\top}$ denotes the transpose of the $3\times3$ matrix $\nabla\widetilde{\d}$.
The vector-valued smooth function $f(\widetilde{\d})$ denotes the penalty function and has the following form:
 $$f(\widetilde{\d})=\nabla_{\widetilde{\d}}F(\widetilde{\d}),$$
where the scalar function $F(\widetilde{\d})$ is the bulk part of the elastic energy.
A typical example is to choose $F(\widetilde{\d})$ as the Ginzburg-Landau penalization thus yielding the penalty function $f(\widetilde{\d})$ as:
$$F(\widetilde{\d})=\frac{1}{4\sigma_0^{2}}(|\widetilde{\d}|^{2}-1)^{2}, \quad  f(\widetilde{\d})=\frac{1}{2\sigma_0^{2}}(|\widetilde{\d}|^{2}-1)\widetilde{\d},$$
where $\sigma_0>0$ is a constant.  We refer the readers to \cite{Chan,DE,E,HK,Leslie1,Lin2} for more physical background and discussion of liquid crystals and mathematical models.

Recall that the Mach number for the compressible flow \eqref{I1} is defined as:
$$M=\frac{|\widetilde{\u}|}{\sqrt{\widetilde{P}'(\widetilde{\rho})}}.$$
 Thus, letting $M$ approach to zero, we hope that $\widetilde{\rho}$, $\widetilde{\d}$ keep a typical size $1$, $\widetilde{\u}$ of the order $\varepsilon$, where $\varepsilon\in (0,1)$ is a small parameter.
We scale $\widetilde{\rho}$, $\widetilde{\u},$ and $\widetilde{\d}$ in the following way:
$$\widetilde{\rho}=\rho_{\varepsilon}(\varepsilon t,x),\quad \widetilde{\u}=\varepsilon\u_{\varepsilon}(\varepsilon t,x),\quad \widetilde{\d}=\d_{\varepsilon}(\varepsilon t,x),$$
and we take the viscosity coefficients as:
$$\widetilde{\mu}=\varepsilon\mu_{\varepsilon},\quad
\widetilde{\lambda}={\varepsilon}^{2}\lambda_{\varepsilon},\quad
\widetilde{\theta}=\varepsilon\theta_{\varepsilon},$$
 where  the normalized coefficients $\mu_{\varepsilon}$, $\lambda_{\varepsilon}$, and $\theta_{\varepsilon}$ satisfy
$$\mu_{\varepsilon}\to \mu,\quad \lambda_{\varepsilon}\to \lambda,\quad
 \theta_{\varepsilon}\to \theta \quad \text{as }\, \varepsilon\to 0^{+},$$
with $\mu$, $\lambda$ and $\theta$ positive constants.
Under this  scaling, system \eqref{I1} becomes
\begin{subequations}\label{I1'}
\begin{align}
&\frac{\partial\rho_{\varepsilon}}{\partial t}+\Dv(\rho_{\varepsilon}\u_{\varepsilon})=0,\label{I1'a}\\
&\frac{\partial(\rho_{\varepsilon} \u_{\varepsilon})}{\partial t}+\Dv(\rho_{\varepsilon} \u_{\varepsilon} \otimes \u_{\varepsilon})+\nabla \frac{1}{\varepsilon^{2}}\rho_{\varepsilon}^{\gamma}\notag\\
&\qquad\qquad\qquad\qquad\qquad= \mu_{\varepsilon} \Delta\u_{\varepsilon}-\lambda_{\varepsilon} \Dv\left(\nabla \d_{\varepsilon} \odot \nabla \d_{\varepsilon}-\big(\frac{1}{2} |\nabla\d_{\varepsilon}|^{2}+F(\d_{\varepsilon})\big) I_{3}\right),\label{I1'b}\\
&\frac{\partial\d_{\varepsilon}}{\partial t}+\u_{\varepsilon}\cdot\nabla\d_{\varepsilon}=\theta_{\varepsilon}(\Delta\d_{\varepsilon}-f(\d_{\varepsilon})),\label{I1'c}
\end{align}
\end{subequations}
where we take $a=1$ because the exact value of $a$ does not play a role in our paper.
The existence of global weak solutions to \eqref{I1'} in bounded domains was established in \cite{WY, LLQ}.
By the initial energy bound \eqref{I7} below, we can assume that the initial datum $\rho_{\varepsilon}^{0}$ is of the order $1+O(\varepsilon)$, so it is reasonable to expect that, as $\varepsilon\to 0$, $\rho_{\varepsilon}\to1$ and  \eqref{I1'a} yields the limit $\Dv\u=0,$ which is the incompressible condition of a fluid, and the first two terms in \eqref{I1'b} become
$$\u_{t}+\Dv(\u\otimes \u)=\u_{t}+(\u\cdot\nabla)\u.$$
The corresponding incompressible equations of liquid crystals  are:
\begin{subequations}\label{I2}
\begin{align}
&\u_{t}+\u\cdot\nabla\u+\nabla \pi= \mu \Delta\u-\lambda \Dv(\nabla \d \odot \nabla \d),\label{I2a}\\
&\d_{t}+\u\cdot\nabla\d=\theta(\Delta\d-f(\d)),\label{I2b}\\
&\Dv\u=0.\label{I2c}
\end{align}
\end{subequations}
Thus, roughly speaking, it is also reasonable to expect from the mathematical point of view that the weak solutions to \eqref{I1'}  converge in  suitable functional spaces to the weak solutions of \eqref{I2} as $\varepsilon\to 0$, and the hydrostatic pressure $\pi$ in \eqref{I2a} is
   the ``limit" of
   $$\frac1{\varepsilon^{2}}(\rho_{\varepsilon}^{\gamma}-1)-\frac{\lambda_{\varepsilon}}2|\nabla\d_{\varepsilon}|^{2}-
   \lambda_{\varepsilon}F(\d_{\varepsilon})$$
    in \eqref{I1'b}. 
This paper is devoted to the rigorous justification of the convergence of
the above incompressible limit (i.e., the low Mach number limit)
for global weak solutions of the compressible  equations of liquid crystals
in  bounded smooth domains.
We remark that the existence of global weak solutions to the incompressible flow of liquid crystals \eqref{I2} was established in Lin-Liu \cite{LL}.

When the direction field $\d$ does not appear, \eqref{I1'} reduces to the compressible Navier-Stokes equations. 
 Lions-Masmoudi \cite{LM} investigated the incompressible limits of the compressible isentropic Navier-Stokes equations in the whole space  and   periodic domains using the group method generated by the wave operator, a method introduced in earlier works \cite{G,SCH} which requires certain smoothness of solutions.
The study in  bounded smooth domains with the no-slip boundary condition on the velocity is much harder
than that in the whole space or periodic domains, because in bounded domains, there are
extra difficulties arising from the appearance of the boundary layers, and the subtle interactions between dissipative effects and wave propagation near the boundary, and hence requires a different approach.
Desjardins-Grenier-Lions-Masmoudi in  \cite{DGLM} relied on spectral analysis and Duhamel's principle to treat these difficulties to find the limit of global solutions in a bounded domain.
These results have been extended by others; see for examples \cite{BDGL,DT, DG,  M, WJ}.
We also remark that in Hoff \cite{Ho} some convergence results were proved for well-prepared data as long as the solution of incompressible limit is suitably smooth. For the case of nonisentropic flows, see \cite{FN,FNP} for some recent developments.  Recently, Hu-Wang \cite{HW} studied the convergence of weak solutions of the compressible magnetohydrodynamic equations to the weak solutions of the incompressible  magnetohydrodynamic equations when 
the Mach number goes to zero in a periodic domain, the whole space, or a bounded domain; and Jiang-Ju-Li \cite{JJL} studied the incompressible MHD limit in the inviscid case as long as the strong solution of the incompressible inviscid MHD exists  with periodic boundary conditions.
For other related studies on the incompressible limits of viscous and inviscid flows,
see \cite{Al, Dan, Ebin, HL, Hoff2, KM1, KM2, Lin, MS1, MS2, Sch1, Sch2}
and the references in \cite{FN}. Finally, we remark that the incompressible flow  can also be derived from the vanishing Debye length type limit of a compressible flow with a
Poisson damping; see for examples \cite{mm1,mm2}.

In this paper, we  shall establish the incompressible limit of \eqref{I1'}
in  a sufficiently smooth bounded domain $\O\subset\R^3$.
As mentioned earlier,  this limit problem in a bounded smooth domain has more difficulties
and requires a different approach due to  the appearance of the boundary
layers and the subtle interactions between dissipative effects and
wave propagation near the boundary.
Comparing with those works on the compressible Navier-Stokes equations, we will encounter extra difficulties in studying the compressible liquid crystals. More precisely, besides the difficulties from compressible Navier-Stokes equations, the appearance of the direction field and the coupling effect between the hodrodynamic equations and the direction field should also be taken into account with new estimates.  We will overcome all these difficulties by adapting the spectral analysis of the semigroup generated by the dissipative wave operator, Duhamel's principle, and the weak convergence method to establish the convergence of the global weak solutions of compressible flow of liquid crystals \eqref{I1'} to the weak solutions of the incompressible flow of liquid crystals \eqref{I2} as 
$\varepsilon$ goes to zero in a bounded domain.

We organize the rest of the paper as follows. In Section 2, we will
provide some preliminaries and state our main result. In Section 3,  we will prove in four steps the convergence of the incompressible limit in a bounded domain.
\bigskip

\section{Preliminaries and Main Results}

We consider the incompressible limit in a smooth bounded domain $\O\subset\R^3$.
To state precisely our main result, we need to introduce a geometrical condition on $\O$ (cf. \cite{DGLM}).
Let us consider the following overdetermined problem:
\begin{equation}\label{I3}-\D\varphi=\lambda\varphi \quad \text{ in } \,\O, \quad \frac{\partial\varphi}{\partial \nu}=0 \quad \text{ on } \,\partial\O, \quad \text{ and } \varphi \text{ is constant on } \partial \O.
\end{equation}
A solution to \eqref{I3} is said to be trivial if $\lambda=0$ and $\varphi$ is a constant. We say that $\O$ satisfies the assumption (H) if all solutions of \eqref{I3} are trivial.
In the two-dimensional case, it was proved that every bounded simply connected open set   with Lipschitz boundary satisfies (H).  We refer the readers to \cite{DGLM} for more information about assumption (H).

Let us  recall the definition of Leray's projectors:  $P$ onto the space of divergence-free vector fields and $Q$ onto the space of gradients, defined by
\begin{equation}
\label{1+}
\u=P\u+Q\u,\quad\quad \text{with} \quad\quad \Dv(P\u)=0,\quad\quad \text{curl}(Q\u)=0,
\end{equation}
for  $\u\in L^{2}.$ Indeed, in view of the results in \cite{Ga}, we know that the operators $P$ and $Q$ are linear bounded operators in $W^{s,p}$ for  $s\geq0$ and $1<p<\infty$ in any bounded domain with smooth boundary.

We consider a sequence of weak solutions $\{(\rho_{\varepsilon},\u_{\varepsilon},\d_{\varepsilon})\}_{\varepsilon>0}$ to \eqref{I1'} in a smooth bounded domain $\O\subset\R^3$ with  the following boundary condition:
\begin{equation}\label{I4}\u_{\varepsilon}|_{\partial\O}=0,\quad\d_{\varepsilon}|_{\partial\O}=\d_{\varepsilon}^{0},\end{equation}
and initial condition:
\begin{equation}\label{I5}\rho_{\varepsilon}|_{t=0}=\rho_{\varepsilon}^{0},\quad \rho_{\varepsilon}\u_{\varepsilon}|_{t=0}=\m_{\varepsilon}^{0},
\quad \d_{\varepsilon}|_{t=0}=\d_{\varepsilon}^{0},\end{equation}
satisfying
\begin{gather}
\rho_{\varepsilon}^{0}\geq 0,\quad \rho_{\varepsilon}^{0}\in L^{\gamma}(\O),\\
\m_{\varepsilon}^{0}\in L^{\frac{2\gamma}{\gamma+1}}(\O),\quad \m_{\varepsilon}^{0}=0 \text{ if }
\rho_{\varepsilon}^{0}=0,\\
\rho_{\varepsilon}^{0}|\u_{\varepsilon}^{0}|^{2}\in L^{1}(\O),\quad \d_{\varepsilon}^{0} \in H^{1}(\O),\\
\left(\rho_{\varepsilon}^{0}\right)^\frac12\u_{\varepsilon}^{0}  \text{ converges weakly in $L^2$ to some $\u_0$  as
$\varepsilon\to 0$,}\\
\d_{\varepsilon}^{0} \text{ converges weakly in $L^2$ to some  $\d_0$  as $\varepsilon\to 0$,}
\end{gather}
and
\begin{equation}\label{I7}
\begin{split}
&\int_{\O}\left(\frac{1}{2}\rho_{\varepsilon}^{0}|\u_{\varepsilon}^{0}|^{2}
+\frac{1}{2}\lambda_{\varepsilon}|\nabla\d_{\varepsilon}^{0}|^{2}+\lambda_{\varepsilon}F(\d_{\varepsilon}^{0})\right)dx
\\&\quad+\frac{1}{\varepsilon^{2}(\gamma-1)}\int_{\O}\left((\rho_{\varepsilon}^{0})^{\gamma}-\gamma\rho_{\varepsilon}^{0}+(\gamma-1)\right)dx \leq C,
\end{split}\end{equation}
for some constant $C>0$.
We remark that \eqref{I7} implies, roughly speaking, that $\rho_{\varepsilon}^{0}$ is of order $1+O(\varepsilon)$
since
$$(\rho_{\varepsilon}^{0})^{\gamma}-\gamma\rho_{\varepsilon}^{0}+(\gamma-1)
=(\rho_{\varepsilon}^{0})^{\gamma}-1-\gamma\left(\rho_{\varepsilon}^{0}-1\right),$$
and $\rho^\gamma$ is a convex function for $\gamma>1$.
As proved in \cite{WY, LLQ},  for any fixed $\varepsilon>0$, there exists a global weak solution $(\rho_{\varepsilon},\u_{\varepsilon},\d_{\varepsilon})$ to the compressible flow of liquid crystals \eqref{I1'} satisfying
$$\rho_{\varepsilon}\in L^{\infty}([0,T];L^{\gamma}(\O)),$$
$$\sqrt{\rho_{\varepsilon}}\u_{\varepsilon} \in L^{\infty}([0,T];L^{2}(\O)),$$
$$\u_{\varepsilon}\in L^{2}([0,T];H^{1}(\O)),$$
$$\d_{\varepsilon} \in L^{2}([0,T];H^{2}(\O))\cap L^{\infty}([0,T];H^{1}(\O)),$$
for any given $T>0$;  and in addition,
$$\rho_{\varepsilon}\u_{\varepsilon}\in C([0,T];L^{\frac{2\gamma}{\gamma+1}}(\O)),$$
$$\rho_{\varepsilon}\in C([0,T];L^{p}_{loc}(\O)),$$
if $1\leq p<\gamma$;  as well as
\begin{equation}\label{I6}E_{\varepsilon}(t)+\mu_{\varepsilon}\int_{0}^{T}\!\!\int_{\O}|\nabla\u_{\varepsilon}|^{2}dxdt
+\lambda_{\varepsilon}\theta_{\varepsilon}\int_{0}^{T}\!\!\int_{\O}|\D\d_{\varepsilon}-f(\d_{\varepsilon})|^{2}dxdt\leq E_{\varepsilon}(0),\end{equation}
for $t \in [0,T]$ a.e., where
$$E_{\varepsilon}:=\int_{\O}\left(\frac{1}{2}\rho_{\varepsilon}|\u_{\varepsilon}|^{2}
+\frac{1}{\varepsilon^{2}(\gamma-1)}\rho_{\varepsilon}^{\gamma}
+\frac{1}{2}\lambda_{\varepsilon}|\nabla\d_{\varepsilon}|^{2}+\lambda_{\varepsilon}F(\d_{\varepsilon})\right)dx,$$
and
$$E_{\varepsilon}(0):=\int_{\O}\left(\frac{1}{2}\rho_{\varepsilon}^{0}|\u_{\varepsilon}^{0}|^{2}
+\frac{1}{\varepsilon^{2}(\gamma-1)}(\rho_{\varepsilon}^{0})^{\gamma}
+\frac{1}{2}\lambda_{\varepsilon}|\nabla\d_{\varepsilon}^{0}|^{2}+\lambda_{\varepsilon}F(\d_{\varepsilon}^{0})\right)dx.$$

We now recall the existence result of global weak solutions to the incompressible flow of liquid crystals in \cite{LL}:
\begin{Proposition}
 \label{p1}For $\u_{0}\in L^{2}(\O)$ and  $\d_{0} \in H^{1}(\O)$ with $\d_{0}|_{\partial\O} \in H^{3/2}(\partial\O)$, system \eqref{I2} with the following initial and boundary conditions:
$$\u|_{t=0}=\u_{0}(x) \quad \text{with }\quad \Dv \u_{0}=0,\quad \d|_{t=0}=\d_{0}(x),$$
and $$\u|_{\partial \O}=0,\quad \d|_{\partial \O}=\d_{0}(x)$$
has a global weak solution $(\u,\d)$
such that
$$\u \in L^{2}(0,T;H^{1}(\O))\cap L^{\infty}(0,T;L^{2}(\O)),$$
$$\d\in L^{2}(0,T;H^{2}(\O))\cap L^{\infty}(0,T;H^{1}(\O)),$$
and
\begin{equation*}
\begin{split}
&\frac{1}{2}\int_{\O}\left(|\u|^2+\lambda|\nabla\d|^2+2\lambda F(\d)\right)dx
+\int_0^{T}\!\!\int_{\O}\left(\mu|\nabla\u|^2+\lambda\theta|\D\d-f(\d)|^2\right)dxdt
\\&\leq\frac{1}{2} \int_{\O}\left(|\u_0|^2+\lambda|\nabla\d_0|^2+2\lambda F(\d_0)\right)dx
\end{split}
\end{equation*}
for all $T \in (0,\infty).$
\end{Proposition}

\smallskip

\begin{Remark}
The global weak solutions obtained in Proposition \ref{p1} are the weak solutions in  Leray's sense. The uniqueness of such solutions can be reached in two-dimensional spaces, see for example \cite{LL}. For more details on the existence and regularity of weak solutions to the incompressible flow of  liquid crystals, we refer the readers to \cite{LL,LL2,HKL}.
\end{Remark}

Our main result reads as follows:

\begin{Theorem}\label{T}
Assume that $\{(\rho_{\varepsilon},\u_{\varepsilon},\d_{\varepsilon})\}_{\varepsilon>0}$ is a sequence of weak solutions to the compressible flow of liquid crystals \eqref{I1'} in a smooth bounded domain $\O\subset\R^3$ with the initial  and boundary conditions \eqref{I4}-\eqref{I7} and $\gamma>\frac{3}{2}.$
Then, for any given $T>0$, as $\varepsilon\to 0$,  $\{(\rho_{\varepsilon},\u_{\varepsilon},\d_{\varepsilon})\}$  converges to a weak solution $(\u,\d)$
of the incompressible flow of liquid crystals \eqref{I2} with the initial data:  $\u|_{t=0}=P\u_{0}$, $\d|_{t=0}=\d_{0}$ and the boundary condition: $\u|_{\partial\O}=0$, $\d|_{\partial\O}=\d_{0}$.  More precisely,
as $\varepsilon\to 0,$
$$\rho_{\varepsilon} \text{ converges to } 1 \text{ in } C([0,T];L^{\gamma}(\O));$$
$$\u_{\varepsilon} \text{ converges to } \u \text { weakly in } L^{2}((0,T)\times\O) \text { and strongly if } \O \text { satisfies condition } (H);$$
$$\d_{\varepsilon} \text { converges to } \d \text{ strongly in } L^{2}([0,T];H^{1}(\O)) \text { and weakly in } L^{2}([0,T];H^{2}(\O)).$$
\end{Theorem}

\bigskip


\section{Proof of Theorem \ref{T}}

In this section we prove Theorem \ref{T} in four steps.

\subsection{A priori estimates and consequences}

We first recall the spectral analysis of the semigroup \cite{DGLM} generated by the dissipative wave operator.
Let $\{\lambda_{k,0}^{2}\}_{k\in \mathbb{N}}\quad(\lambda_{k,0}>0)$ be the nondecreasing sequence of eigenvalues and  $\{\Phi_{k,0}\}_{k\in \mathbb{N}}$ in $L^{2}(\O)$ be 
the eigenvectors with zero mean value  of the Laplace operator satisfying the  homogeneous Neumann boundary condition:
$$-\D\Phi_{k,0}=\lambda_{k,0}^{2}\Phi_{k,0}\; \text{ in } \O,  \quad \frac{\partial \Phi_{k,0}}{\partial \nu}=0 \; \text { on } \partial \O,$$
where $\nu$ is the unit outer normal of $\O$.
By the Gram-Schmidt orthogonalization method, it is possible to assume that $\{\Phi_{k,0}\}_{k\in\mathbb{N}}$ is an orthonormal basis of $L^{2}(\O)$ and that up to a slight modification, if $\lambda_{k,0}=\lambda_{l,0},$ and $k\neq l$, then $$\int_{\partial \O}\nabla \Phi_{k,0}\cdot\nabla\Phi_{l,0} ds=0.$$

From \eqref{I6} and the conservation of mass, we have for almost all $t\geq0$,
\begin{equation}\label{p1}\begin{split}&\int_{\O}\left(\frac{1}{2}\rho_{\varepsilon}|\u_{\varepsilon}|^{2}+\frac{1}{\varepsilon^{2}(\gamma-1)}(\rho_{\varepsilon}^{\gamma}-\gamma
\rho_{\varepsilon}+\gamma-1)+\frac{\lambda_{\varepsilon}}{2}|\nabla_{\varepsilon}\d_{\varepsilon}|^{2}+\lambda F(\d_{\varepsilon})\right)dx
\\& \quad
+\mu_{\varepsilon}\int_{0}^{T}\!\!\int_{\O}|\nabla\u_{\varepsilon}|^{2}dxdt
+\lambda_{\varepsilon}\theta_{\varepsilon}\int_{0}^{T}\!\!\int_{\O}|\D\d_{\varepsilon}-f(\d_{\varepsilon})|^{2}dxdt\\
&\leq \int_{\O}\left(\frac{1}{2}\rho_{\varepsilon}^{0}|\u_{\varepsilon}^{0}|^{2}+\frac{1}{\varepsilon^{2}(\gamma-1)}((\rho_{\varepsilon}^{0})^{\gamma}-\gamma
\rho_{\varepsilon}^{0}+\gamma-1)+\frac{\lambda_{\varepsilon}}{2}|\nabla_{\varepsilon}\d_{\varepsilon}^{0}|^{2}+\lambda F(\d_{\varepsilon}^{0})\right)dx\\&\leq C.
\end{split}\end{equation}
By \eqref{p1}, we have the following properties:
\begin{equation}
\label{p1+}
\sqrt{\rho_{\varepsilon}}\u_{\varepsilon} \text{ is bounded in } L^{\infty}([0,T];L^{2}(\O)),
\end{equation}
\begin{equation}
\label{p2+}
\d_{\varepsilon} \text { is bounded in } L^{\infty}([0,T];H^{1}(\O)),
\end{equation}
\begin{equation}
\label{p3+}
\frac{1}{\varepsilon^{2}(\gamma-1)}(\rho_{\varepsilon}^{\gamma}-\gamma\rho_{\varepsilon}
+\gamma-1) \text { is bounded in } L^{\infty}([0,T];L^{1}(\O));
\end{equation}
and
\begin{equation}
\label{p4+}
\nabla\u_{\varepsilon} \text { is bounded in } L^{2}([0,T];L^{2}(\O)),
\end{equation}
\begin{equation}
\label{p5+}
\D\d_{\varepsilon}-f(\d_{\varepsilon}) \text { is bounded in } L^{2}([0,T];L^{2}(\O))
\end{equation}
for all $T>0$.

Let us recall the following basic fact:
\begin{subequations}
\text { for some } $c_0>0$ \text{ and for all } $x\geq 0$,
\begin{align}
&x^{\gamma}-1-\gamma(x-1)\geq c_0|x-1|^{2} \text{ if } \gamma\geq2,\notag\\
&x^{\gamma}-1-\gamma(x-1)\geq c_0|x-1|^{2} \text{ if } \gamma<2 \text { and } x\leq R,\notag\\
& x^{\gamma}-1-\gamma(x-1)\geq c_0|x-1|^{\gamma} \text{ if }\gamma <2 \text { and }x\geq R,\notag
\end{align}
\end{subequations}
where $R\in (0,\infty).$
Thus from \eqref{p3+}, we have
\begin{equation}
\label{p6+}
\int_{\O}\left(\frac{1}{\varepsilon^{2}}|\rho_{\varepsilon}-1|^{2}\chi_{|\rho_{\varepsilon}-1|\leq1/2}
+\frac{1}{\varepsilon^{2}}|\rho_{\varepsilon}-1|^{\gamma}\chi_{|\rho_{\varepsilon}-1|\geq1/2}\right)dx \leq C,
\end{equation}
where  $ \chi$ is the characteristics function and $C$ denotes a generic positive constant hereafter.
By \eqref{p6+}, one has
\begin{equation*}
\sup_{t\geq0}\|\rho_{\varepsilon}-1\|_{L^{\gamma}(\O)}\leq C\varepsilon^{\kappa/\gamma},\quad\quad \text { and } \quad\quad\quad \sup_{t\geq0}\|\rho_{\varepsilon}-1\|_{L^{\kappa}(\O)}\leq C \varepsilon,
\end{equation*}
where $\kappa=\min\{2,\gamma\}$,
which implies that
 \begin{equation}
 \label{p10+}
\rho_{\varepsilon} \to 1 \text { in } C([0,T];L^{\gamma}(\O))\quad\quad  \text{ as } \varepsilon \to 0.
\end{equation}

Now we split the velocity
$$\u_{\varepsilon}=\u_{\varepsilon}^{1}+\u_{\varepsilon}^{2}\quad \text { with } \quad \u_{\varepsilon}^{1}=\u_{\varepsilon}\chi_{|\rho_{\varepsilon}-1|\leq \frac{1}{2}}, \quad \u_{\varepsilon}^{2}=\u_{\varepsilon}\chi_{|\rho_{\varepsilon}-1|>\frac{1}{2}},$$
which means that
\begin{equation}
\label{p8+}
\sup_{t\geq0}\int_{\O}|\u_{\varepsilon}^{1}|^{2}dx\leq 2\sup_{t\geq0}\int_{\O}\rho_{\varepsilon}|\u_{\varepsilon}|^{2}dx \leq C
\end{equation}
and \begin{equation}
\label{p9+}
\|\u_{\varepsilon}^{2}\|_{L^{2}(\O)}^{2}\leq 2\int_{\O}|\rho_{\varepsilon}-1||\u_{\varepsilon}|^{2}dx\leq C\varepsilon \|\u_{\varepsilon}\|_{L^{2\kappa/(\kappa-1)}(\O)}^{2}\leq C\varepsilon\|\nabla\u_{\varepsilon}\|_{L^{2}(\O)}^{2},
\end{equation}
where we used some embedding inequality.
By \eqref{p8+} and \eqref{p9+}, it is easy to see that
\begin{equation*}
\u_{\varepsilon}^{1} \text{ is bounded in } L^{\infty}([0,T];L^{2}(\O)),
 \end{equation*}
and
\begin{equation*}
 \u_{\varepsilon}^{2}\varepsilon^{-1/2} \text{ is bounded in } L^{2}([0,T];L^{2}(\O)).
  \end{equation*}
This implies that $\u_{\varepsilon}$ is bounded in $L^{2}((0,T)\times\O)$ for all $T>0.$

By \eqref{I1'a} and \eqref{p10+}, one deduces that,
\begin{equation*} 
\Dv \u_{\varepsilon}\to 0 \text { weakly in } L^{2}([0,T];L^{2}(\O))
\end{equation*}
for all $T>0,$
where we used a fact
\begin{equation*}
\int_{0}^{T}\!\!\int_{\O}|\Dv\u_{\varepsilon}|^{2}dxdt \leq C\int_{0}^{T}\!\!\int_{\O}|\nabla\u_{\varepsilon}|^{2}dxdt \leq C.
\end{equation*}
By \eqref{p5+}, smoothness of $f$, and the standard elliptic theory, we have
\begin{equation*}
\nabla^{2}\d_{\varepsilon} \in L^{2}([0,T];L^{2}(\O)).
\end{equation*}
On the other hand, multiplying $\d_{\varepsilon}$ on the both sides of \eqref{I1'c}, then applying maximal principle, we have
 $\d_{\varepsilon} \in L^{\infty}([0,T]\times\O).$
 Thus, \begin{equation*}
 \d_{\varepsilon} \in L^{2}([0,T];H^{2}(\O)).
 \end{equation*}
Similarly to \cite{WY}, using the Gagliardo-Nirenberg inequality,
$$\|\nabla\d_{\varepsilon}\|_{L^{4}(\O)}\leq c\|\Delta\d_{\varepsilon}\|_{L^{2}(\O)}^{\frac{1}{2}}\|\d_{\varepsilon}\|_{L^{\infty}(\O)}^{\frac{1}{2}}+c\|\d_{\varepsilon}\|_{L^{\infty}(\O)},$$
one has
\begin{equation}\label{p2}\nabla\d_{\varepsilon} \in L^{4}((0,T)\times\O).\end{equation}
Summing up the above estimates, we can assume that, up to a subsequence if necessary,
$$\u_{\varepsilon}\to \u \text { weakly in } L^{2}([0,T];H^{1}(\O),$$
$$\Dv\u_{\varepsilon}\to 0 \text { weakly in }L^{2}([0,T];L^{2}(\O)), $$
$$\rho_{\varepsilon}\to 1 \text{ in }C([0,T];L^{\gamma}(\O)),$$
 $$\D\d_{\varepsilon}-f(\d_{\varepsilon})\to \D\d-f(\d) \text{ weakly in } L^{2}([0,T];L^{2}(\O)),$$
and $$\d_{\varepsilon}\to \d \text { weakly in } L^{\infty}(0,T;H^{1}(\O))\cap L^{2}([0,T];H^{2}(\O)).$$
\smallskip

To show the strong convergence of $\d_{\varepsilon}$, we rely on the following Aubin-Lions compactness lemma (see \cite{Lj}):

\begin{Lemma}\label{L}
Let $X_{0}, X$ and $X_{1}$ be three Banach spaces with $X_{0}\subseteq X\subseteq X_{1}$
and $X_{0}$ and $X_{1}$ be reflexive. Suppose that $X_{0}$ is compactly embedded in $X$ and that $X$ is continuously embedded  in $X_{1}$. For $1 <p,q<\infty,$ let
$$W=\left\{u\in L^{p}([0,T];X_{0}):\; \frac{du}{dt} \in L^{q}([0,T];X_{1})\right\}.$$
Then the embedding of $W$ into $L^{p}([0,T];X)$ is also compact.
 \end{Lemma}

From \eqref{I1'c}, it is easy to see that
\begin{equation}
\begin{split}
&\|\partial_{t}\d_{\varepsilon}\|_{L^{2}(\O)}
\leq\|\theta_{\varepsilon}(\D\d_{\varepsilon}-f(\d_{\varepsilon}))\|_{L^{2}(\O)}
+\|\u_{\varepsilon}\cdot\nabla\d_{\varepsilon}\|_{L^{2}(\O)}\\& \leq C\|\u_{\varepsilon}\|^{2}_{L^{4}(\O)}+C\|\nabla\d_{\varepsilon}\|^{2}_{L^{4}(\O)}+C\|\D\d_{\varepsilon}-f(\d_{\varepsilon})\|_{L^{2}(\O)},
\\& \leq C\|\nabla \u_{\varepsilon}\|^{2}_{L^{2}(\O)}+C\|\nabla\d_{\varepsilon}\|^{2}_{L^{4}(\O)}+C\|\D\d_{\varepsilon}-f(\d_{\varepsilon})\|_{L^{2}(\O)}.
\end{split}\end{equation}
This, combined with \eqref{p4+}, \eqref{p5+}, \eqref{p2} and embedding theorem, implies that
$$\|\partial_{t}\d_{\varepsilon}\|_{L^{2}([0,T];L^{2}(\O))}\leq C.$$
Since $H^{2}\subset H^{1}\subset L^{2}$ and the injection $H^{2}\hookrightarrow H^{1}$ is compact, we apply Lemma \ref{L} to deduce that the sequence $\{\d_{\varepsilon}\}$ is precompact in $L^{2}(0,T;H^{1}(\O)).$
 By taking a subsequence if necessary, we can assume that,
 $$\d_{\varepsilon}\to\d \text { weakly in } L^{2}(0,T;H^{2}(\O)),$$
and  $$\d_{\varepsilon}\to\d \text { strongly in } L^{2}(0,T;H^{1}(\O)).$$
Therefore, by a standard argument, we deduce that the limit $\u$ and $\d$ satisfy equation \eqref{I2c} in the sense of distributions.
 By  strong convergence of $\d_{\varepsilon}$ and smoothness of $F$, we deduce the convergence of the nonlinear term
 \begin{equation*}
 \begin{split}
 &\nabla\d_{\varepsilon}\odot\nabla\d_{\varepsilon}-(\frac{1}{2}|\nabla\d_{\varepsilon}|^{2}+F(\d_{\varepsilon}))I \\
 &\to
 \nabla\d\odot\nabla\d-(\frac{1}{2}|\nabla\d|^{2}+F(\d))I \quad \quad \text{ as } \varepsilon\to 0
 \end{split}
 \end{equation*}
 in the sense of distributions.
Denoting
$$\pi_{\varepsilon}=\frac{1}{\varepsilon^{2}}\rho_{\varepsilon}^{\gamma}-\frac{\lambda_{\varepsilon}}{2}|\nabla\d_{\varepsilon}|^{2}-\lambda_{\varepsilon}F(\d_{\varepsilon}),$$
 then we rewrite   equation \eqref{I1'b} as
\begin{equation}\label{p3}\frac{\partial(\rho_{\varepsilon} \u_{\varepsilon})}{\partial t}+\Dv(\rho_{\varepsilon} \u_{\varepsilon} \otimes \u_{\varepsilon})+\nabla \pi_{\varepsilon}= \mu_{\varepsilon} \Delta\u_{\varepsilon}-\lambda_{\varepsilon} \Dv(\nabla \d_{\varepsilon} \odot \nabla \d_{\varepsilon}).\end{equation}
 We project equation \eqref{p3} onto divergence-free vector fields:
 \begin{equation}\label{p4}\partial_{t}P(\rho_{\varepsilon}\u_{\varepsilon})+P(\Dv(\rho_{\varepsilon}\u_{\varepsilon}\otimes\u_{\varepsilon}))
 -\mu_{\varepsilon}\D P\u_{\varepsilon}=-P(\lambda_{\varepsilon}\Dv(\nabla\d_{\varepsilon}\odot\nabla\d_{\varepsilon})),\end{equation}
where $P$ is defined by \eqref{1+}.
 Then \eqref{p4} yields a bound on $\partial_{t}P(\rho_{\varepsilon}\u_{\varepsilon})$ in
 \begin{equation*}
 L^{2}([0,T];H^{-1}(\O))+L^{2}([0,T];W^{-1,1}(\O))+L^{1}([0,T];H^{-1}(\O)),
 \end{equation*}
and hence in $L^{1}([0,T];W^{-1,1}(\O)).$
In addition, $P(\rho_{\varepsilon}\u_{\varepsilon})$ is bounded in
$$L^{\infty}([0,T];L^{\frac{2\gamma}{\gamma+1}}(\O))\cap L^{2}([0,T];L^{r}(\O)) \text{ with }
\frac{1}{r}=\frac{1}{\gamma}+\frac16.$$

 To continue our proof, we need the following lemma (cf. Lemma 5.1 in \cite{L}).

 \begin{Lemma}\label{l2}
 Let the functions $ g_{n}$, $h_{n} $ converge weakly to the functions $g$, $h$, respectively, in $L^{p_{1}}(0,T;L^{p_{2}}(\O)),$ $L^{q_{1}}(0,T;L^{q_{2}}(\O))$, where $1\leq p_{1},p_{2}\leq\infty$ and
 $$\frac{1}{p_{1}}+\frac{1}{q_{1}}=\frac{1}{p_{2}}+\frac{1}{q_{2}}=1.$$
 Assume, in addition, that
 $$\frac{\partial g_{n}}{\partial t} \text { is bounded in } L^{1}(0,T;W^{-m,1}(\O)) \text { for some } m\geq0 \text { independent of }n,$$
  and $$\|h_{n}-h_{n}(t,\cdot+\zeta)\|_{L^{q_{1}}(0,T;L^{q_{2}})}\to 0 \text { as } |\zeta|\to 0 \text { uniformly in } n.$$
  Then $g_{n}h_{n}$ converges to $gh$ in the sense of distributions in $(0,T)\times \O.$
  \end{Lemma}

  Applying Lemma \ref{l2}, we deduce that $P(\rho_{\varepsilon}\u_{\varepsilon})\cdot P(\u_{\varepsilon})$ converges to $|\u|^{2}$ in the sense of distributions.
   It is easy to see that $P\u_{\varepsilon} \to \u=P\u \text{ in the sense of distributions}$
  because the weak convergence of $P\u_{\varepsilon}$ to $\u=P\u$ in $L^{2}((0,T);L^{2}(\O))$
   and $$\int_{0}^{T}\!\!\int_{\O}(|P\u_{\varepsilon}|^{2}-P(\rho_{\varepsilon}\u_{\varepsilon})\cdot P\u_{\varepsilon})dxdt \leq C\|\rho_{\varepsilon}-1\|_{C([0,T];L^{\gamma}(\O))}\|\u_{\varepsilon}\|_{L^{2}([0,T];L^{s}(\O))}^{2},$$
   with $s=\frac{2\gamma}{\gamma-1}<6$ since $\gamma>\frac{3}{2}$.

\subsection{The convergence of $Q\u_{\varepsilon}$}
To prove our main result, it remains to show the convergence of the gradient part of the velocity $Q\u_{\varepsilon}.$  The argument for the convergence of $Q\u_{\varepsilon}$ in our paper follows the same line in \cite{DGLM} and \cite{HW}, except  the argument for the direction field $\d_{\varepsilon}$. For the convenience of readers and the completeness of argument, we provide the details here.

First, we introduce the spectral problem associated with the viscous wave operator $L_{\varepsilon}$ in terms of eigenvalues and eigenvectors of the invisicid wave operator $L$. In the sequel, we write the density fluctuation as
$$\varphi_{\varepsilon}=\frac{\rho_{\varepsilon}-1}{\varepsilon},$$
and $\phi=(\Phi,\m)^{\top} $.
We define the wave operators $L$ and $L_{\varepsilon}$ in $\mathcal D'(\O)\times \mathcal D'(\O)$ as follows:
\begin{equation}
L\left(\begin{array}{c}
\Phi\\ \m\end{array}\right)=\left(\begin{array}{c}\Dv \m \\ \nabla \Phi
\end{array} \right)
\end{equation}
and \begin{equation} L_{\varepsilon}\left(\begin{array}{c}
\Phi \\ \m\end{array}\right)= L\left(\begin{array}{c}\Phi\\ \m\end{array}\right)+\varepsilon \left(\begin{array}{c} 0 \\ \mu_{\varepsilon}\D \m
\end{array}\right).\end{equation}
The eigenvalues and eigenvectors of $L$ read as follows:
\begin{equation}
\phi_{k,0}^{\pm}=\left(\begin{array}{c}\Phi_{k,0}\\ \m_{k,0}^{\pm}=\pm\frac{\nabla\Phi_{k,0}}{i\lambda_{k,0}}
\end{array}\right),\end{equation}
and $$L\phi_{k,0}^{\pm}=\pm i\lambda_{k,0}\phi_{k,0}^{\pm} \text{ in } \O,\quad \m_{k,0}^{\pm}\cdot\nu=0 \text {  on } \partial\O.$$
And we need the following lemma \cite{DGLM}:

\begin{Lemma}\label{l3}
Let $\O$ be a $C^{2}$ bounded domain of $\R^{d}$ and let $k\geq1,N\geq0$. Then, there exist approximate eigenvalues $i\lambda_{k,\varepsilon,N}^{\pm}$ and eigenvectors $\phi_{k,\varepsilon,N}^{\pm}=(\Phi_{k,\varepsilon,N}^{\pm},m_{k,\varepsilon,N}^{\pm})^{\top}$ of $L_{\varepsilon}$ such that $$L_{\varepsilon}\phi_{k,\varepsilon,N}^{\pm}
=i\lambda_{k,\varepsilon,N}^{\pm}\phi_{k,\varepsilon,N}^{\pm}+R_{k,\varepsilon,N}^{\pm},$$
with
$$ i\lambda_{k,\varepsilon,N}^{\pm}=\pm i\lambda_{k,0}+i\lambda_{k,1}^{\pm}\sqrt{\varepsilon}+O(\varepsilon), \text{ where } Re(i\lambda_{k,1}^{\pm})\leq 0, $$
and for all $p \in [1,\infty], $ we have,
$$\left|R_{k,\varepsilon,N}^{\pm}\right|_{L^{p}(\O)}\leq C_{p}(\sqrt{\varepsilon})^{N+1/p}\quad \text { and }\quad \left|\phi_{k,\varepsilon,N}-\phi_{k,0}^{\pm}\right|_{L^{p}(\O)}\leq C_{p}(\sqrt\varepsilon)^{1/p}. $$
\end{Lemma}
\bigskip

\begin{Remark}
\label{l0}
The key idea is to construct an approximation scheme of $L_{\varepsilon}$ in terms of $\phi_{k,0}^{\pm}$. We refer the readers to \cite{DGLM} for more details and the proof. From Lemma \ref{l3} and its proof, we have, for any integers $k$,
\begin{equation*}
i\lambda_{k,1}^{\pm}=-\frac{1\pm i}{2}\sqrt{\frac{\mu_{\varepsilon}}{2\lambda_{k,0}^{3}}}\int_{\partial\O}|\nabla\Phi_{k,0}|^2dx,
\end{equation*}
which satisfies
\begin{equation}
\label{0+}
 Re (i\lambda_{k,1}^{\pm})\leq 0.
\end{equation}
We observe that the first order term $i\lambda_{k,1}^{\pm}$ clearly yields an instantaneous damping of acoustic waves, as soon as $Re (i\lambda_{k,1}^{\pm})< 0.$ Thus, we let $$I\subset\mathbb{N}$$
to be a collection of the all eigenvectors $\Phi_{k,0}$ of the Laplace operator such that
$$Re (i\lambda_{k,1}^{\pm})< 0.$$
 Denote $$J=\mathbb{N}\setminus I,$$ that is to say, when $k\in J$, we have $$Re (i\lambda_{k,1}^{\pm})=0$$ due to \eqref{0+}. This implies $\lambda_{k,0}=0.$ In the case that $\lambda_{k,0}=0$, $\m_{k,0}^{\pm}$ must vanish on $\partial\O$ and therefore not only
$\m_{k,0}^{\pm}\cdot \nu=0$ but also $\m_{k,0}^{\pm}=0$ on $\partial \O.$ Thus, no significant boundary layer is created, and there is no enhanced dissipation of energy in these layers.
\end{Remark}

Remark that $\{\frac{\nabla\Phi_{k,0}}{\lambda_{k,0}}\}_{k \in \mathbb N}$ is an orthonormal basis of $L^{2}(\O)$ functions with zero mean value on $\O$. We write
\begin{equation*}
Q\u_{\varepsilon}=\sum_{k\in \mathbb N}\left(Q \u_{\varepsilon},\frac{\nabla\Phi_{k,0}}{\lambda_{k,0}}\right)\frac{\nabla\Phi_{k,0}}{\lambda_{k,0}}, \end{equation*}
where $$(f,g)=\int_{\O}f(x)\overline{g(x)}dx.$$
We split $Q\u_{\varepsilon}$ into two parts $Q_{1}\u_{\varepsilon}$ and $Q_{2}\u_{\varepsilon}$, defined by
$$Q_{1}\u_{\varepsilon}=\sum_{k\in I } \left(Q\u_{\varepsilon}, \frac{\nabla\Phi_{k,0}}{\lambda_{k,0}}\right)\frac{\nabla\Phi_{k,0}}{\lambda_{k,0}},$$
and $$Q_{2}\u_{\varepsilon}=\sum_{k\in J } \left(Q\u_{\varepsilon}, \frac{\nabla\Phi_{k,0}}{\lambda_{k,0}}\right)\frac{\nabla\Phi_{k,0}}{\lambda_{k,0}},$$
which, respectively, correspond to damped terms and nondamped terms.
What remains is to show that
 \begin{equation*}
Q_{1}\u_{\varepsilon}\to 0 \quad\quad\text { in }L^{2}([0,T]\times\O),
 \end{equation*}
 and  \begin{equation*}
 \text{curl}\,\Dv(Q_{2}\m_{\varepsilon}\otimes Q_{2}\u_{\varepsilon}) \to 0
  \end{equation*}
in the sense of distributions if $J\neq\emptyset,$ which is equivalent to saying that $\Dv(Q_{2}\m_{\varepsilon}\otimes Q_{2}\u_{\varepsilon})$
converges to a gradient in the sense of distributions.

Let us observe that in view of the bound on $\u_{\varepsilon}$ in $L^{2}(0,T;H^{1}_{0}(\O)),$ the problem reduces to a finite number of modes. Remark that the eigenvalues $\{\lambda_{k,0}^{2}\}_{k\geq0}$ is a nondecreasing sequence, we have
\begin{equation*}
\sum_{k>N}\int_{0}^{T}\left|\left(Q_{i}\u_{\varepsilon},\frac{\nabla\Phi_{k,0}}{\lambda_{k,0}}\right)\right|^{2}dt\leq \frac{C}{\lambda_{N+1}^{2}}\left|\nabla\u_{\varepsilon}\right|_{L^{2}((0,T)\times\O)}^{2}, \quad i=1 \text{ or } 2.
\end{equation*}
Letting $ N\to \infty$, then $\lambda_{N}\to \infty$, which implies that $(Q_{1}\u_{\varepsilon},\m_{k,0}^{\pm})\to 0 \text { in } L^{2}(0,T)$ for any $k>N.$  So we need to show that $(Q_{1}\u_{\varepsilon},\m_{k,0}^{\pm})$ converges to $0$ strongly in $ L^{2}(0,T)$ for any fixed $k$ and study the interaction of a finite number of terms in $\Dv(Q_{2}\u_{\varepsilon}\otimes Q_{2}\u_{\varepsilon}).$

Recalling that $\varphi_{\varepsilon}=\frac{\rho_{\varepsilon}-1}{\varepsilon}$, we have $$Q\u_{\varepsilon}=Q\m_{\varepsilon}-\varepsilon Q(\varphi_{\varepsilon}\u_{\varepsilon})$$
and \begin{equation}\begin{split}
\varepsilon|(Q(\varphi_{\varepsilon}\u_{\varepsilon}),\nabla\Phi_{k,0})|
&=\varepsilon\left|\int_{\O}\varphi_{\varepsilon}\u_{\varepsilon}\cdot\nabla\Phi_{k,0}dx\right|
\\& \leq \varepsilon\|\varphi_{\varepsilon}\|_{L^{\gamma}(\O)}\|\u_{\varepsilon}\|_{L^{\frac{\gamma}{\gamma-1}}(\O)}\|\nabla\Phi_{k,0}\|_{L^{\infty}(\O)},
\end{split}\end{equation}
which tends to zero in $L^{2}(0,T)$ due to $\gamma>\frac{3}{2}.$ So we move to study $(Q\m_{\varepsilon},\m_{k,0}^{\pm}).$
We write $$\beta_{k,\varepsilon}^{\pm}=(\phi_{\varepsilon}(t),\phi_{k,0}^{\pm})$$with
\begin{equation}
\phi_{\varepsilon}(t)=\left(\begin{array}{c}\varphi_{\varepsilon}\\ \m_{\varepsilon}
\end{array}\right).
\end{equation}
It is easy to see $$2(Q\m_{\varepsilon},\m_{k,0}^{\pm})=\beta_{k,\varepsilon}^{\pm}-\beta_{k,\varepsilon}^{\mp}.$$
Applying Lemma \ref{l3} with $N=2$ and the H\"{o}lder inequality, one obtains
\begin{equation}\begin{split}
&|(\phi_{\varepsilon}(t),\phi_{k,0}^{\pm})-\phi_{k,\varepsilon,2}^{\pm}|\leq \|\phi_{\varepsilon}\|_{L^{2}(\O)}\|\phi_{k,0}^{\pm}-\phi_{k,\varepsilon,2}^{\pm}\|_{L^{2}(\O)}\\ &
\leq C \varepsilon^{\frac{\alpha}{2}}\left(\|\varphi_{\varepsilon}\|_{L^{\infty}([0,T];L^{\kappa}(\O))}+
\|\m_{\varepsilon}\|_{L^{\infty}([0,T];L^{\frac{2\gamma}{\gamma+1}}(\O))}\right)
\end{split}\end{equation}
where $$\alpha=\min\{1-\frac{1}{\kappa},\frac{1}{2}-\frac{1}{2\gamma}\}.$$
It remains to show that $b_{k,\varepsilon}^{\pm}(t)=(\phi_{\varepsilon}(t),\phi_{k,\varepsilon,2}^{\pm})$
converges to zero strongly in $L^{2}([0,T])$ when $k \in I$ and check the oscillations when $k\in J.$

Using $L_{\varepsilon}^{*}$ to denote the adjoint operator of $L_{\varepsilon}$ with respect to $(\cdot,\cdot)$, we have
\begin{equation}
\label{p5}
\partial_{t}\phi_{\varepsilon}-\frac{L_{\varepsilon}^{*}\phi_{\varepsilon}}{\varepsilon}
=\left(\begin{array}{c}0\\g_{\varepsilon}
\end{array}\right)\end{equation}
where $$g_{\varepsilon}=-\Dv(\m_{\varepsilon}\otimes \m_{\varepsilon})-\nabla\pi_{\varepsilon}-\lambda_{\varepsilon}\Dv(\nabla\d_{\varepsilon}\odot\nabla\d_{\varepsilon}).$$
Taking the scalar product of \eqref{p5} with $\phi_{k,\varepsilon,2}^{\pm}$, one obtains
\begin{equation}\label{p11+}
(\partial_{t}\phi_{\varepsilon},\phi_{k,\varepsilon,2}^{\pm})-
\left(\frac{L_{\varepsilon}^{*}\phi_{\varepsilon}}{\varepsilon},\phi_{k,\varepsilon,2}^{\pm}\right)
=c_{k,\varepsilon}^{\pm}(t)
\end{equation}
where $c_{k,\varepsilon}^{\pm}(t)=(g_{\varepsilon},\m_{k,\varepsilon,2}^{\pm})+\varepsilon^{-1}(\phi_{\varepsilon},R_{k,\varepsilon,2}^{\pm}).$

Letting $b_{k,\varepsilon}^{\pm}(t)=(\phi_{\varepsilon},\phi_{k,\varepsilon,2}^{\pm})$ and observing that
$$\left(\frac{L_{\varepsilon}^{*}\phi_{\varepsilon}}{\varepsilon},\phi_{k,\varepsilon,2}^{\pm}\right)
=\frac{1}{\varepsilon}\overline{i\lambda_{k,\varepsilon}^{\pm}}(\phi_{\varepsilon},\phi_{k,\varepsilon,2}^{\pm})
+\frac{1}{\varepsilon}R_{k,\varepsilon,2}^{\pm},$$
we rewrite \eqref{p11+} as
\begin{equation}\label{p6}
\partial_{t}b_{k,\varepsilon}^{\pm}-
\frac{1}{\varepsilon}\overline{i\lambda_{k,\varepsilon,2}^{\pm}}b_{k,\varepsilon}^{\pm}=c_{k,\varepsilon}^{\pm}.
\end{equation}
We can view \eqref{p6} as an ordinary differential equation, which satisfies the initial value:
\begin{equation}
\label{p12+}
b_{k,\varepsilon}^{\pm}|_{t=0}=b_{k,\varepsilon}^{\pm}(0).
\end{equation}

\subsection{The case $k \in I$.}

By the basic theory of ordinary differential equation, the solution to \eqref{p6}-\eqref{p12+} is given by
\begin{equation}\label{p7}b_{k,\varepsilon}^{\pm}(t)
=b_{k,\varepsilon}^{\pm}(0)e^{\overline{i\lambda_{k,\varepsilon,2}^{\pm}}\frac{t}{\varepsilon}}
+\int_{0}^{t}c_{k,\varepsilon}^{\pm}(s)e^{\overline{i\lambda_{k,\varepsilon,2}^{\pm}}\frac{(t-s)}{\varepsilon}}ds.\end{equation}
Next, we need to estimate $b_{k,\varepsilon}^{\pm}(t)$.
Note that $$i\lambda_{k,\varepsilon,2}^{\pm}=\pm i\lambda_{k,0}+i\lambda_{k,1}^{\pm}\sqrt{\varepsilon}+o(\varepsilon),$$
which is helpful in estimating the first part of \eqref{p7}:
\begin{equation}\|b_{k,\varepsilon}^{\pm}(0)e^{\overline{i\lambda_{k,\varepsilon,2}}\frac{t}{\varepsilon}}\|_{L^{2}(0,T)}
\leq C\|b_{k,\varepsilon}^{\pm}(0) e^{Re(\overline{i\lambda_{k,1}^{\pm}})\frac{t}{\sqrt{\varepsilon}}}\|_{L^{2}(0,T)}
\leq C\varepsilon^{\frac{1}{4}},\end{equation}
where we used Lemma \ref{l3}.

To estimate the second term in \eqref{p7}, observe that for any $a\in L^{q}(0,T)$ and $1\leq p,q\leq\infty \text { such that } \frac{1}{p}+\frac{1}{q}=1,$ we have
\begin{equation}\label{p8}
\left|\int_{0}^{T}e^{\overline{i\lambda_{k,\varepsilon,2}}\frac{t-s}{\varepsilon}}a(s)ds\right|
\leq\int_{0}^{T}\left| e^{Re(\overline{i\lambda_{k,1}^{\pm}})\frac{t-s}{\sqrt{\varepsilon}}}\right||a(s)|ds\leq C|a|_{L^{q}(0,T)}\sqrt{\varepsilon}^{\frac{1}{p}}.
\end{equation}
To prove
\begin{equation*}
b_{k,\varepsilon}^{\pm}(t)\to 0 \text{ strongly in } L^2([0,T]),
\end{equation*}
it remains to show that $c_{k,\varepsilon}^{\pm}$ is bounded in $L^{q}(0,T)$ for some $q>1.$
Remark that $$|c_{k,\varepsilon}^{\pm}|\leq c_{1}+c_{2}+c_{3}+c_{4},$$
where $$c_{1}=\left|\int_{\O}(\m_{\varepsilon}\otimes \m_{\varepsilon})(t)\cdot\nabla \m_{k,\varepsilon,2}^{\pm}dx\right|$$
and $c_2, c_3, c_4$ will be defined as follows.
Using $$\m_{\varepsilon}=\varepsilon\varphi_{\varepsilon}\u_{\varepsilon}+\u_{\varepsilon},\quad
\u_{\varepsilon}=\u_{\varepsilon}^{1}+\u_{\varepsilon}^{2},$$  and integration by parts,
we have
\begin{equation}
\begin{split}\label{p13+}
& c_{1}(t)\leq \varepsilon \int_{\O}|\varphi_{\varepsilon}\u_{\varepsilon}\otimes\u_{\varepsilon}\cdot\nabla \m_{k,\varepsilon,2}^{\pm}|dx+C\int_{\O}|\u_{\varepsilon}\cdot\nabla\u_{\varepsilon}\cdot \m_{k,\varepsilon,2}^{\pm}|dx
\\& \leq \|\m_{k,\varepsilon,2}^{\pm}\|_{L^{\infty}(\O)}\|\u_{\varepsilon}^{1}+\u_{\varepsilon}^{2}\|_{L^{2}(\O)}\|\nabla\u_{\varepsilon}\|_{L^{2}(\O)}
\\ &\quad +\varepsilon\|\varphi_{\varepsilon}\|_{L^{\infty}([0,T];L^{\kappa}(\O))}\|(\u_{\varepsilon})^{2}\|_{L^{\kappa/\kappa-1}(\O)}\|\|\nabla\m_{k,\varepsilon,2}^{\pm}\|_{L^{\infty}(\O)}\\
& \leq C\|\u_{\varepsilon}^{1}\|_{L^{\infty}([0,T];L^{2}(\O))}\cdot\|\nabla\u_{\varepsilon}\|_{L^{2}(\O)}
+C\varepsilon^{1/2}\|\nabla\u_{\varepsilon}\|_{L^{2}(\O)}^{2}+C\varepsilon^{1/2}\|\nabla\u_{\varepsilon}\|_{L^{2}(\O)}.
 \end{split}
 \end{equation}
For the second term
\begin{equation*}
c_{2}=\int_{\O}(\pi_{\varepsilon}\Dv \m_{k,\varepsilon,2}^{\pm})dx,
\end{equation*}
recalling that
\begin{equation*}
\Dv \m_{k,\varepsilon,2}^{\pm}=i \lambda_{k,\varepsilon,2}^{\pm}\Phi_{k,\varepsilon,2}+R_{k,\varepsilon,2}^{\pm},
\end{equation*}
 we have \begin{equation}
\label{P14+}
c_{2}(t)\leq\|\pi_{\varepsilon}\|_{L^{\infty}([0,T];L^{1}(\O))}(\|\Phi_{k,\varepsilon,2}^{\pm}\|_{L^{\infty}(\O)}
+\|R_{k,\varepsilon,2}^{\pm}\|_{L^{\infty}(\O)})\leq C.
\end{equation}
The third term can be estimated as:
\begin{equation} \label{p15+}
\begin{split}
&c_{3}(t)=\varepsilon^{-1}|(\phi_{\varepsilon},R_{k,\varepsilon,2}^{\pm})|\\
&\leq \varepsilon^{-1}\|R_{k,\varepsilon,2}^{\pm}\|_{L^{\kappa/(\kappa-1)}(\O)}\|\phi_{\varepsilon}\|_{L^{\infty}([0,T];L^{\kappa}(\O))}\\
&\leq C\varepsilon^{1/2-1/2\kappa}.
\end{split}
\end{equation}
And the last term is
\begin{equation}
\begin{split}
\label{p16+}
&c_{4}(t)=\left|\int_{\O}\Dv(\nabla\d_{\varepsilon}\odot\nabla\d_{\varepsilon})\m_{k,\varepsilon,2}^{\pm}dx\right|
\\ &\leq C\int_{\O}|\nabla\d_{\varepsilon}|^{2} |\nabla \m_{k,\varepsilon,2}^{\pm}|dx
\\ &\leq C\|\nabla\m_{k,\varepsilon,2}^{\pm}\|_{L^{\infty}(\O)}\|\nabla\d_{\varepsilon}\|_{L^{2}(\O)}
\\ &\leq C \|\nabla\d_{\varepsilon}\|_{L^{2}(\O)},
\end{split}
\end{equation}
where we used  $\nabla\d_{\varepsilon} \in L^{4}(\O).$
From \eqref{p13+}-\eqref{p16+}, it is easy to see that the term $c_{k,\varepsilon}^{\pm}$ is bounded, and consequently, $$b_{k,\varepsilon}^{\pm}\to 0 \text { strongly in } L^{2}(0,T).$$
To complete our proof, we need to consider further the case $k\in J$.

\subsection{ The case $k\in J$.}
As in Remark \ref{l0}, when $k\in J$, we have $\lambda_{k,1}^{\pm}=0.$
This implies that, together with \eqref{p7},
$$e^{\pm i\lambda_{k,0}t/\varepsilon}b_{k,\varepsilon}^{\pm}\; \text { is bounded in } L^{2}(0,T),$$
 and $$\partial_{t}(e^{\pm i\lambda_{k,0}t/\varepsilon}b_{k,\varepsilon}^{\pm}) \text { is bounded in } \sqrt{\varepsilon}L^{1}(0,T)+ L^{p}(0,T) \text { for some } p>1.$$
 It follows that, up to a subsequence if necessary, $b_{k,\varepsilon}^{\pm}$ converges strongly in $L^{2}(0,T)$ to some element $b_{k,osc}^{\pm}$.
 Since $\rho_{\varepsilon}\to 1 \text { in } C([0,T];L^{\gamma}(\O))$ and $b_{k,\varepsilon}^{\pm}$ is uniformly bounded in $L^{2}(0,T),$
one obtains that $$\rho_{\varepsilon}b_{k,\varepsilon}^{\pm}b_{l,\varepsilon}^{\pm}\frac{\nabla\Phi_{k,0}}{\lambda_{k,0}}\otimes\frac{\nabla\Phi_{l,0}}{\lambda_{l,0}}
-b_{k,\varepsilon}^{\pm}b_{l,\varepsilon}^{\pm}\frac{\nabla\Phi_{k,0}}{\lambda_{k,0}}\otimes\frac{\nabla\Phi_{l,0}}{\lambda_{l,0}}\to 0 $$ in the sense of distributions. Thus, we only need to consider the terms
$$b_{k,\varepsilon}^{\pm}b_{l,\varepsilon}^{\pm}\frac{\nabla\Phi_{k,0}}{\lambda_{k,0}}\otimes\frac{\nabla\Phi_{l,0}}{\lambda_{l,0}} \quad\quad \text { for all } k,l \in J.$$
On the other hand,
\begin{equation*}\begin{split}& b_{k,\varepsilon}^{\pm}b_{l,\varepsilon}^{\pm}\frac{\nabla\Phi_{k,0}}{\lambda_{k,0}}\otimes\frac{\nabla\Phi_{l,0}}{\lambda_{l,0}}
\\ & =e^{i(\lambda_{k,0}-\lambda_{l,0})t/\varepsilon}e^{-i\lambda_{k,0}t/\varepsilon}b_{k,0}e^{i\lambda_{l,0}t/\varepsilon}b_{l,0}
\frac{\nabla\Phi_{k,0}}{\lambda_{k,0}}\otimes\frac{\nabla\Phi_{l,0}}{\lambda_{l,0}} \\& \to
e^{i(\lambda_{k,0}-\lambda_{l,0})t/\varepsilon}b_{k,osc}^{\pm}b_{l,osc}^{\pm}\frac{\nabla\Phi_{k,0}}
{\lambda_{k,0}}\otimes\frac{\nabla\Phi_{l,0}}{\lambda_{l,0}} \quad\quad\text {as } \varepsilon \to 0
\end{split}\end{equation*}
in the sense of distributions,
where we used the strong convergence of $$e^{\pm \lambda_{k,0}t/\varepsilon}b_{k,\varepsilon}^{\pm} \text { in } L^{2}(0,T) \quad\quad\quad\text { when } k\in J.$$

If $\lambda_{k,0}=\lambda_{l,0},$
we can write
 $$\Dv(\nabla\Phi_{k,0}\otimes\nabla\Phi_{l,0}+\nabla\Phi_{l,0}\otimes\nabla\Phi_{k,0})
 =-\lambda_{k,0}^{2}\nabla(\Phi_{k,0}\Phi_{l,0})+\nabla(\nabla\Phi_{k,0}\cdot\nabla\Phi_{l,0}),$$
 which is a gradient, and thus disappears in the pressure term. It remains to consider the case $\lambda_{k,0}\neq \lambda_{l,0}$.
  In this case, we have a fact that
  $$ b_{k,osc}^{\pm}(t)b_{l,osc}^{\pm}(t) \in L^{1}([0,T])$$ for all $k,l \in J$ due to
 $$b_{k,osc}^{\pm}(t) \in L^{2}([0,T]) \quad \quad \text { for all } k \in J.$$
 This implies, together with the Riemann-Lebesgue Lemma,
 \begin{equation*}
 \int_{\O}e^{i(\lambda_{k,0}-\lambda_{l,0})t/\varepsilon}b_{k,osc}^{\pm}(t)b_{l,osc}^{\pm}(t)dt\to 0 \quad\quad \text { as } \varepsilon \to 0.
 \end{equation*}
Thus, we conclude that
\begin{equation*}
e^{i(\lambda_{k,0}-\lambda_{l,0})t/\varepsilon}b_{k,osc}^{\pm}b_{l,osc}^{\pm}\frac{\nabla\Phi_{k,0}}
{\lambda_{k,0}}\otimes\frac{\nabla\Phi_{l,0}}{\lambda_{l,0}} \to 0  \quad\quad \text { as } \varepsilon\to 0
\end{equation*}
in the sense of distributions.  Hence, the finite sum, as $k,l \leq N,$
$$\Dv\left(\sum_{k,l\in J}\rho_{\varepsilon}b_{k,\varepsilon}^{\pm}b_{l,\varepsilon}^{\pm}
\frac{\nabla\Phi_{k,0}}{\lambda_{k,0}}\otimes\frac{\nabla\Phi_{l,0}}{\lambda_{l,0}}\right)$$
converges to a gradient in the sense of distributions, which means that
$$\Dv\left(\rho_{\varepsilon}Q_{2}\u_{\varepsilon}\otimes Q_{2}\u_{\varepsilon}\right)$$
converges to a gradient in the sense of distributions.

The proof of Theorem \ref{T} is now complete.

\bigskip\bigskip

\section*{Acknowledgments}

D. Wang's research was supported in part by the National Science
Foundation under Grant DMS-0906160 and by the Office of Naval
Research under Grant N00014-07-1-0668. C. Yu's research was supported in part by the National Science
Foundation under Grant DMS-0906160.

\end{document}